\documentclass[12pt]{article}
\usepackage{amsfonts,amssymb}
\usepackage{latexsym}
\usepackage{theorem}
\usepackage{amsmath}
\newcommand{\linespacing}{1.25}
\renewcommand{\baselinestretch}{\linespacing}
\newtheorem{theorem}{Theorem}[section]
\newtheorem{lemma}[theorem]{Lemma}
\newtheorem{proposition}[theorem]{Proposition}
\newtheorem{corollary}[theorem]{Corollary}

{\theorembodyfont{\normalfont\rmfamily}
\newtheorem{definition}[theorem]{Definition}

\newtheorem{example}{Example}
\newtheorem{remark}[theorem]{Remark}
}

\makeatletter
\def\ack{\vspace{.5\baselineskip}\noindent{\theorem@headerfont
Acknowledgement}\ \ } %%
\newenvironment{introthm}[1]%
    {\renewcommand{\theththm}{\ref{#1}}\begin{ththm}}
    {\end{ththm}}
\newtheorem{ththm}{Theorem}
\newenvironment{proof}[1][]%
{\def\proof@temp{#1}\par\noindent
\textsc{Proof}\ifx\proof@temp\@empty\else\
({\proof@temp})\fi\hspace{1em}} {\hphantom{xxx}\hfill~
{$\Box$}\par\vspace{.4\baselineskip}} %%
\def\operatorname#1{\mathop{\operator@font #1}\nolimits}%
\makeatother %%
\newcommand{\map}[1]{\stackrel{#1}{\longrightarrow}}
\newcommand{\R}{\mathbb{R}}
\renewcommand{\L}{\mathcal{L}}
\newcommand{\G}{\mathcal{G}}
\newcommand{\g}{\mathfrak{g}}
\newcommand{\h}{\mathfrak{h}}
\renewcommand{\H}{\mathcal{H}}

\newcommand{\Tr}{\operatorname{Tr}}
\newcommand{\Ker}{\operatorname{Ker}\,}

\newcommand{\half}{{\textstyle{\frac12}}}
\def\cyclic{\mathop{\kern0.9ex{{+}
\kern-2.2ex\raise-.28ex\hbox{\Large\hbox
{$\circlearrowright$}}}}\limits}

\def\ftnote#1{\def\footnotemark{}\footnote{#1}\setcounter{footnote}{0}}

\begin{document}
  \title{Reduction,Induction and Ricci flat
symplectic connections.}
 \author{Michel Cahen${}^{*}$
 and 
 Simone Gutt${}^{* , **}$
\ftnote{\kern-4.5pt${}^{*}$Universit\'e Libre de Bruxelles, Campus
Plaine CP 218, Bvd du Triomphe, B-1050~Brussels, Belgium}
\ftnote{\kern-4.5pt${}^{**}$Universit\'e de Metz, D\'epartement de math\'ematiques, Ile du Saulcy,
F-57045~Metz Cedex 01, France} \ftnote{Email:
\texttt{mcahen@ulb.ac.be}, \texttt{sgutt@ulb.ac.be}}   }
\date{ We are pleased to dedicate this paper to Hideki Omori on the
occasion of his $65^{th}$ birthday}
\setcounter{page}{0}
 %% Don't want front page spaced out
 \renewcommand{\baselinestretch}{1}
 \maketitle
\thispagestyle{empty}
\begin{abstract}
 In this paper we
present a construction of Ricci-flat connections  through an induction
procedure. Given a symplectic manifold $(M,\omega)$ of dimension $2n$, 
we define induction as a way
to construct a symplectic manifold $(P,\mu)$ of dimension $2n+2$. Given
any symplectic connection $\nabla$ on $(M,\omega)$, we define an induced
connection $\nabla^P$ which is a Ricci-flat symplectic connection on $(P,\mu)$.
 \end{abstract}
  \newpage
\renewcommand{\baselinestretch}{\linespacing}
 %%%%%%%%%%%%%%%%
 %%%%%%%%%%%%%%%%%
 \section*{Introduction}\label{sect:intro}
 
 A symplectic connection on a symplectic manifold $(M,\omega)$
 is a torsionless linear connection $\nabla$ 
 on $M$ for which the symplectic $2$--form $\omega$ is parrallel.
 A symplectic connection exists on any symplectic manifold and the space 
 of such connections is an
 affine space modelled on the space of symmetric $3$--tensorfields on $M$.
 
 In all what follows, the dimension $2n$ of the manifold $M$ is
 assumed to be $\ge 4$ unless explicitely stated.
 The curvature tensor $R^\nabla$ of a symplectic connection $\nabla$ decomposes
 \cite{Vaisman} under the action of the symplectic group into $2$ irreducible
 components, $R^\nabla=E^\nabla+W^\nabla$.
The $E^\nabla$ component is defined only in terms of the
 Ricci-tensor $r^\nabla$ of $\nabla$. All traces of the $W^\nabla$
component vanish. 

Two particular types of symplectic connections thus arize:

- symplectic connections for which $W^\nabla=0;$ we call them Ricci-type symplectic connections;  

- symplectic connections  for which $E^\nabla=0$; they are
called   Ricci-flat since $E^\nabla=0 \Leftrightarrow  r^\nabla=0$.

When studying \cite{CGS} local and global models for Ricci-type symplectic connections,
(or more generally \cite{CS} so called special 
symplectic connections) ,  Lorenz Schwachh\"{o}fer
and the present authors were lead to  consider examples of the following
construction:\\
$\bullet$ start with a symplectic manifold  $(M,\omega)$ of dimension $2n$;\\
$\bullet$ build a (cooriented) contact manifold $(N,\alpha)$ of dimension $2n+1$ and a
submersion $\pi:N\rightarrow M$ such that  $d\alpha =\pi^*\omega$;\\
$\bullet$ define on the manifold  $P=N\times \mathbb R$ a natural
 symplectic structure $\mu$.\\
It was
observed  \cite{CGS} that    if $(M,\omega)$ admits a  symplectic
connection of Ricci type one could ``lift'' this  connection to $P$
and the lifted connection is symplectic (relative  to $\mu$) and
flat. 

The aim of this paper is to generalize this  result. More
precisely we formalize the notion of induction for symplectic manifolds. 
Starting from a symplectic manifold $(M,\omega)$, we define a 
contact quadruple
$(M, N,\alpha, \pi)$, where $N,\alpha$ and $\pi$ are as above,
 and  we build   the corresponding $2n+2$ dimensional
symplectic manifold $(P,\mu)$.
We prove the following: 
\begin{introthm}{th:Ricciflat}  Let $(M,\omega)$ be a symplectic
manifold which is the first element of  a contact quadruple $(M, N,
\alpha, \pi)$. Let $\nabla$ be an  artitrary symplectic connection on
$(M,\omega)$. Then one can lift  $\nabla$ to a symplectic connection
on $(P,\mu)$ which is Ricci--flat.
\end{introthm}
This theorem has various applications. In particular one has
\begin{introthm}{th:indred}
Let $(P,\mu)$ be a symplectic manifold admitting 
 a conformal vector field $S$ which is
complete,  a symplectic vector field $E$ which  commutes
with $S$ and assume that, for any $x\in P,~\mu_x(S,E)>0$.
Assume the reduction of $\Sigma=\{x\in P\,|\,\mu_x(S,E)=1\}$
by the flow of $E$ has a manifold structure $M$ with $\pi:\Sigma\rightarrow M$
a surjective submersion.
Then $(P,\mu)$ admits a Ricci-flat connection.
 \end{introthm}
  The paper is organized as follows. In
section \ref{sec:1}  we  
 study sufficient conditions 
for a symplectic manifold $(M,\omega)$ to be the first element of a
contact quadruple and we give examples of such quadruples.
Section \ref{sec:2} is devoted to the lift of hamiltonian
 (resp conformal) vector fields from $(M,\omega)$ 
 to the induced symplectic manifold $(P,\mu)$ constructed via a contact
 quadruple.
We  show that if $(M,\omega)$ is conformal homogeneous,
so is $(P,\mu)$. Section \ref{sec:3} describes the structure
of conformal homogeneous symplectic manifolds; this part is
certainly known but as we had no immediate reference we decided to
include it. Section \ref{sec:4} gives some constructions of
lifts of symplectic connections of $(M,\omega)$
to symplectic connections on the induced symplectic manifold $(P,\mu)$ 
constructed via a contact
 quadruple. We also
prove theorem \ref{th:Ricciflat}.
In section \ref{sec:5} we give conditions for a symplectic manifold
$(P,\mu)$ to be obtained by induction from a contact quadruple $(M,N,\alpha,\pi)$.
 We  give also a proof of theorem \ref{th:indred}. 

\section{Induction and contact quadruples}\label{sec:1}

  \begin{definition}
 A {\bf{contact quadruple}} is a quadruple $(M, N,\alpha, \pi)$
 where $M$ is a $2n$ dimensional smooth  manifold,
$N$ is a smooth $2n+1$ dimensional manifold, $\alpha$ is a cooriented
 contact structure on $N$ (i.e. $\alpha$ is a $1$--form on $N$ such that
 $\alpha\wedge (d\alpha)^n$ is nowhere vanishing), $\pi: N\rightarrow M$
 is a smooth submersion and $d\alpha =\pi^*\omega$ where $\omega$
 is a symplectic $2$--form on $M$.
 \end{definition}
 \begin{definition}
 Given a contact quadruple $(M, N,\alpha, \pi)$
  the {\bf{induced symplectic manifold}} is the $2n+2$ dimensional manifold
  $$P:=N\times \R$$ endowed with the (exact) symplectic structure
 $$
\mu:=2e^{2s}~ds\wedge p_1^*\alpha +e^{2s}~ dp_1^*\alpha =d(e^{2s}~p_1^*\alpha)
$$
 where $s$ denotes the variable along $\R$ and  
$p_1: P \rightarrow N$ the projection on the first factor. 
\end{definition}
Induction in the sense of building a $2n+2$-- dimensional symplectic manifold
from a symplectic manifold of dimension $2n$ is also considered
by Kostant in \cite{Kostant}.\\
\begin{remark}
$\bullet$ The vector field $S:=\partial_s$ on $P$ is such that
$i(S)\mu=2{\rm e}^{2s}(p_1^*\alpha)$; hence 
$L_{S}\mu=2\mu$ and $S$ is a conformal vector
field.\\
$\bullet$ The Reeb vector field $Z$ on $N$ (i.e. the vector field $Z$ on $N$
such that $i(Z)d\alpha=0$ and $i(Z)\alpha=1$) lifts to a vector field $E$
 on $P$ such that:
$p_{1*}E=Z$ and $ds(E)=0.$
Since $i(E)\mu=-d({\rm e}^{2s})$, $E$ is a Hamiltonian
vector field on $(P,\mu)$.  Furthermore
\begin{eqnarray*}
   [E,S] &=& 0 \\
   \mu(E,S) &=& -2{\rm e}^{2s}.
\end{eqnarray*}
$\bullet$ Observe also that if $\Sigma=\{\, y\in P\, |\, s(y)=0\, \}$, the reduction of
 $(P,\mu$) relative to the constraint manifold $\Sigma$ (which is
 isomorphic to $N$) is precisely $(M,\omega)$.\\
$\bullet$ For $y\in P$ define $H_y(\subset T_yP)= >E,S<^{\bot_\mu}$.
Then $H_y$ is symplectic and $(\pi\circ p_1)_{*y}$ defines a linear
isomorphism between $H_y$ and $T_{\pi p_1(y)}M$. Vector fields on
$M$ thus admit
 ``horizontal'' lifts to $P$.\\
 \end{remark}
We shall now make some remarks on the existence of a contact quadruple 
the first term of which corresponds to a given symplectic manifold $(M,\omega)$.

\begin{lemma}\label{lem1.1}  Let $(M,\omega)$ be a smooth
symplectic manifold of dimension $2n$  and let $N$ be a smooth
$(2n+1)$ dimensional manifold admitting a  smooth surjective
submersion $\pi$ on $M$. Let ${\mathcal H}$ be a  smooth $2n$
dimensional distribution on $N$ such that $\pi_{*x}:  {\mathcal
H}_x\rightarrow T_{\pi(x)}M$ is a linear isomorphism (remark that
such a distribution   may always be constructed by
choosing a smooth riemannian metric $g$ on $N$ and setting
${\mathcal H}_x=(\ker \pi_{*x})^\bot$). Then  either
there exists a smooth nowhere vanishing 1-form $\alpha$ and a smooth
vector field $Z$ such that $\forall x\in N$ we have (i) $\ker\alpha_x={\mathcal H}_x$ (ii) $Z_x\in$
$\ker\pi_{*x}$ (iii) $\alpha_x(Z_x)=1$ or the same is true for a  double
cover of $N$.
 \end{lemma}
 \begin{proof} Choose an auxiliary riemannian
metric $g$ on $M$ and consider $N'=\{Z\in TN\,|\,
Z\in\ker\pi_* {\rm{~and~} }g(Z,Z)=1\}$. If $N'$ has two components,
one can choose a global vector field $Z\in \ker\pi_*$ on $N$ and 
define a smooth 1-form $\alpha$ with $\ker\alpha={\mathcal H}$ 
and $\alpha(Z)=1$. 
If $N'$ is connected, $N'$ is a double  cover of $N$
($p:N'\rightarrow N: Z_x\rightarrow x)$ and we can choose  coherently $Z'\in
T_ZN'$ by the rule that its projection on $T_xN$  is precisely $Z$.
\end{proof}
 This says that if we have a pair $(M,N)$ with a surjective
submersion $\pi:N\rightarrow M$ we can always assume (by passing
eventually to a double cover of $N$) that there exists a nowhere
vanishing vector field $Z\in \ker\pi_*$ and a nowhere vanishing
1-form $\alpha$ such that $\alpha(Z)=1$ and $\ker\alpha$ projects
isomorphically on the tangent space to $M$. The vector field $Z$ is
 determined up to non zero multiplicative factor by the submersion
 $\pi$; on the other hand, having chosen $Z$, the 1-form $\alpha$
can  be modified by the addition of an arbitrary 1-form $\beta$
vanishing  on $Z$. If $\tilde\alpha=\alpha+\beta$ is another choice,
the 2-form  $d\tilde\alpha$ is the pull back of a 2-form on $M$ iff
 $i(Z)d\tilde\alpha=0$; i. e. iff:
  $$ \mbox{(i) }
L_Z\beta=-{ L}_Z\alpha\qquad\qquad\mbox{(ii)  }\beta(Z)=0.
$$
  This can always be solved locally. We shall assume this can
be  solved globally.
\begin{lemma}\label{lem1.2}
 Let $(M,\omega)$ be a smooth symplectic manifold of dimension $2n$
 and let $N$ be a smooth $(2n+1)$ dimensional manifold admitting a
  smooth surjective submersion $\pi$ on $M$. Let $Z$ be a smooth
  nowhere vanishing vector field on $N$ belonging to $\ker\pi_*$. Let
  $\alpha$ be a
1-form such that $\alpha(Z)=1$. If ${ L}_Z\alpha=\mu\alpha$, for
a certain $\mu\in C^\infty(N)$,  then $\mu=0$ and $d\alpha$ is
the pullback of a closed 2-form $\nu$  on $M$. Furthermore if $X$
(resp $Y$) is a vector field on $M$ and  $\bar X$ (resp $\bar Y$) is
the vector field on $N$ such that (i)  $\pi_*\bar X=X$ (resp.
$\pi_*\bar Y=Y$) (ii) $\alpha(\bar  X)=\alpha(\bar Y)=0$ then:
\begin{eqnarray*}
   [\bar X,\bar Y]-\overline{[X,Y]} &=& -\pi^*(\nu(X,Y))Z \\
   {[}Z,\bar X] &=& 0.
\end{eqnarray*}
\end{lemma}
\begin{proof} We have $\pi_*[ Z,\bar X]=0$, $[Z,\bar X]=\alpha([Z,\bar X])Z=
-({ L}_Z\alpha)(\bar X)Z=0$.\\
Since $({ L}_Z\alpha)(Z)=d\alpha(Z,Z)=0$,
$\mu$ vanishes. Also:
$$
i(Z)d\alpha={\cal L}_Z\alpha=0.
$$
so $d\alpha$ is the pullback of a closed $2$--form $\nu$ on $M$. Finally:
 \begin{eqnarray*}
   \pi_*[\bar X,\bar Y]&=&\pi_ *\overline{[X,Y]} \\
   {[}\bar X,\bar Y]&=& \overline{[X,Y]}+\alpha([\bar X,\bar
Y])Z=\overline{[X,Y]}-d\alpha(\bar X,\bar Y)Z.
\end{eqnarray*}
\end{proof}
\begin{corollary}\label{corol1.2}
 If $\nu=\omega$, the manifold
$(N,\alpha)$ is a contact manifold and  $Z$ is the corresponding
Reeb vector.
\end{corollary}
We shall now give examples of contact quadruples for given symplectic manifolds.
\begin{example}
Let $(M, \omega=d\lambda)$ be an exact symplectic  manifold.
Define $N=M\times\mathbb R$, $\pi=p_1$
(=projection of the  first factor), $\alpha=dt+p_1^*\lambda$;
then $(N,\alpha)$ is a contact  manifold and $(M,N,\alpha,\pi)$
is a contact quadruple.

The associated induced manifold is
$P=N\times\mathbb R=M\times\mathbb R^2$; with coordinates 
$(t,s)$ on $\R^2$ and obvious identification
$$  \mu={\rm e}^{2s}~[\,d\lambda +
2ds\wedge  (dt+\lambda)\,].
$$
\end{example}

\begin{example} Let $(M,\omega)$
be a quantizable symplectic manifold; this means that there is a complex
line bundle $L\map{p}M$ with hermitean structure $h$ and a connection $\nabla$ 
on $L$ preserving $h$ whose
curvature is proportional to $i\omega$.\\
Define $N: = \{\,\xi\in L~|~h(\xi,\xi\,)=1\}\subset L$ to be the
unit circle sub-bundle. It is a principal  $U(1)$ bundle and $L$ is
the associated bundle $L=N\times_{U(1)}\mathbb C$. The connection
1-form on $N$ (representing $\nabla$) is $u(1)=i\R$ valued
and will be denoted $\alpha'$; its curvature
is $d\alpha'=ik\omega$. Define $\alpha:=\dfrac{1}{ik}\alpha'$ 
and $\pi:=p|_N : N\rightarrow M$  the surjective submersion.
Then $(M,N,\alpha,\pi)$ is a contact quadruple.

The associated induced manifold $P$ is in bijection with
$L_0=L\setminus$ zero section; indeed, consider
 $$
\Psi:L_0\rightarrow P=N\times\mathbb R:
\xi\rightarrow\left(\frac{\xi}{h(\xi,\xi)^{1/2}},k\ln
h(\xi,\xi)^{1/2}\right).
$$
Clearly $L_0$ is a $\mathbb C^*$
principal bundle on $M$; denote by $\check{\alpha}$ the $\mathbb C^*$--valued
$1$--form on $L_0$ representing $\nabla$; if $j_1:  N\rightarrow L_0$ is the natural
injection and similarly $j_2:  i\mathbb R\rightarrow \mathbb C$ the
obvious injection, we have
$$
j_1^*\check{\alpha}=j_2\circ \alpha'.
$$
Then
\begin{eqnarray*}
\left((\Psi^{-1})^*\check{\alpha}\right)_{(\xi_0,s)}(X_{\xi_0}+a\partial_s)
&=& \check{\alpha}_{\xi_0{\rm e}^{s/k}}(\Psi^{-1}_*(X_{\xi_0}+a\partial_s)  \\
&=& \check{\alpha}_{\xi_0}({R_{{\rm
e}^{-s/k}}}_*\circ\Psi^{-1}_*(X_{\xi_0}+a\partial_s))\\
&=&\check{\alpha}_{\xi_0}(X_{\xi_0}+a\partial_s)=j_2\, \alpha'(X_{\xi_0})+\frac{a}{k}
\end{eqnarray*}
i. e.
$$
\Psi^{-1*}\check{\alpha}=p_1^*j_2\, \alpha'+\frac{ds}{k}.
$$
On the other hand the 1-form e$^{2s}p_1^*\alpha=\dfrac{1}{ik}{\rm
e}^{2s}p_1^*j_2^*\alpha'\,$; this shows how the symplectic
form $\mu=d({\rm e}^{2s}~p_1^*\alpha)$ 
on $P$ is related to the connection form on $L_0$ 
[$\mu=d(\frac{e^{2s}}{ik}\Psi^{-1*}\check{\alpha}$)].\\
Such examples have been studied by Kostant \cite{Kostant}.
\end{example}

\begin{example} Let $(M,\omega)$ be a connected homogeneous symplectic
 manifold; i. e. $M=G/H$ where $G$ is a Lie group which we may
assume  connected and simply connected and where $H$ is the
stabilizer in  $G$ of a point $x_0\in M$. If $p: G\rightarrow
M:g\rightarrow gx_0,\quad p^*\omega$ is a left invariant closed 2-form
on $G$ and  $\Omega=(p^*\omega)_e$, (e=neutral element of $G$) is a
Chevalley  2-cocycle on ${\g}$ (=Lie Algebra of $G$) with values
in $\mathbb R$ (for the trivial representation).\\
Notice that $\Omega$ vanishes as soon as one  of its arguments is
 in ${\h}$ (=Lie algebra of $H$). Let ${\g}_1=\g
\oplus\mathbb R$ be the central extension of ${\g}$ defined by
 $\Omega$; i. e.
 $$
  [(X,a),(Y,b)]=([X,Y],\Omega (X,Y)).
  $$
 Let ${\h}'$ be the subalgebra of $\g_1$, isomorphic to $\h$,  
defined by ${\h}':=\{\,(X,0)\,\vert\,X\in \h\,\}$. Let $G_1$
be the connected and simply connected group of algebra $\g$, 
and let $H'$ be the connected subgroup of $G_1$
with Lie algebra ${\h}'$. {\bf{Assume $H'$ is closed}}.\\
Then $G_1/H'$ admits a natural structure of smooth manifold;
define $N:=G_1/H'$.
Let $p_1:G_1\rightarrow G$
be the homomorphism whose differential is the  projection
$\g_1\rightarrow \g$ on the first factor; clearly $p_1(H')\subset H$.
Define $\pi:N=G_1/H'\rightarrow M=G/H : g_1H'\mapsto p_1(g_1)H$; it is a surjective
submersion.

We shall now construct the contact form $\alpha$ on $N$:
  $p_1^*\circ p^*\omega$ is a left
invariant  closed 2-form on $G_1$ vanishing on the fibers of $p\circ
p_1:G_1\rightarrow M$. Its  value $\Omega_1$ at the neutral element $e_1$ of
$G_1$ is a  Chevalley 2-cocycle of $\g_1$ with values in $\mathbb R$.
Define the $1$--cochain  $\alpha_1: \g_1\rightarrow \mathbb R: (X,a)\rightarrow
-a$. Then
\begin{eqnarray*}
\Omega_1((X,a),(Y,b))&=&
(p^*\omega)_{e}(X,Y)
=\Omega(X,Y)=-\alpha_1([(X,a),(Y,b)])\\
&=&\delta\alpha_1((X,a),(Y,b)),
\end{eqnarray*}
i. e. $\Omega_1=\delta\alpha_1$ is a coboundary. Let
$\tilde\alpha_1$ be the left invariant 1-form on $G_1$ corresponding
 to $\alpha_1$. Let $q: G_1\rightarrow G_1/H'=N$ be the natural projection.
 We shall show that there exists a $1$--form $\alpha$ on $N$
 so that $q^*\alpha=\tilde\alpha_1$.\\
  For any $U\in \g_1$ denote by $\tilde U$  the
corresponding  left invariant vector field on $G_1$.
 For any $X\in \h'$ we have
\begin{eqnarray*}
   &&i(\tilde X)\tilde\alpha_1 = \alpha_1(X)=0 \\
   &&(L_{\tilde X}\tilde\alpha_1)(\widetilde (Y,b)) =
-\tilde\alpha_1([\tilde X,\widetilde (Y,b)])=-\alpha_1([X,(Y,b)])=\Omega(X,Y)=0 
\end{eqnarray*}
so that indeed $\tilde\alpha_1$ is the pullback by $q$ of a $1$--form $\alpha$ on $N=G_1/H'$.
Furthermore $d\alpha=\pi^*\omega$ because both are $G_1$ invariant $2$--forms on $N$ and:
\begin{eqnarray*}
 (d\alpha)_{q(e_1)}((X,a)^{*N},(Y,b)^{*N})&=
 &(q^*d\alpha)_{e_1}(\widetilde{(X,a)},\widetilde{(Y,b)})\\
 &=&(d\tilde\alpha_1)_{e_1}(\widetilde{(X,a)},\widetilde{(Y,b)})\\
 &=&\Omega(X,Y),\\
 &=&\omega_{x_0}(X^{*M},Y^{*M})\\
 &=&(\pi^*\omega)_{q(e_1)}((X,a)^{*N},(Y,b)^{*N})
  \end{eqnarray*}
where we denote  by $U^{*N}$  the 
 fundamental vector field on $N$ associated to $U \in \g_1$  . 
 Thus
\begin{lemma}\label{lem1.3}
Let ($M=G/H,\omega$) be a homogeneous symplectic manifold; let
$\Omega$ be the value at the neutral element of $G$ of the pull back
 of $\omega$ to $G$. This is a Chevalley 2 cocycle of the Lie algebra ${\g}$ 
 of $G$. If ${\g}_1=\g\oplus \R$ is the central extension
of  ${\g}$ defined by this 2 cocycle and $G_1$ is the
corresponding  connected and simply connected group let $H'$ be the
connected  subgroup of $G_1$ with algebra ${\h}'=\{\,(X,0)\,\vert\,X\in \h\,\}\cong\h$.
Assume $H'$ is a closed subgroup of $G_1$. Then $N=G_1/H'$ admits a
natural submersion $\pi$ on $M$ and has a contact structure $\alpha$
such that $d\alpha=\pi^*\omega$. 
Hence $(G/H,G_1/H',\alpha,\pi)$ is a contact quadruple.
\end{lemma}
\begin{remark} The
center of $G_1$ is connected and simply connected, hence the central subgroup
exp$t(0,1)$ is isomorphic to $\mathbb R$. The subgroup $p_1^{-1}(H)$ is a closed
Lie
 subgroup of $G_1$ whose connected component is $p_1^{-1}(H_0)$
($H_0$ =connected component of $H$).   The universal cover
$\widetilde{p_1^{-1}(H_0)}$ of $p_1^{-1}(H_0))$ is
 the  direct product of $\tilde H_0$ (=universal cover of
$H_0$)  by $\mathbb R$. If $\nu:\widetilde{p_1^{-1}(H_0)}\rightarrow
 p_1^{-1}(H_0)$ is the covering homomorphism, the subgroup $H'$
 we are interested in is $H'=\nu(\tilde H_0)$.
 Clearly if $\pi_1(H_0)\sim \ker \nu$ is finite , $H'$ is closed and 
 the construction proceeds.
 \end{remark}
\end{example}

\section{Lift of hamiltonian vector fields and of conformal vector
 fields}\label{sec:2}

Let $(M,\omega)$ be a symplectic manifold of dimension $2n$ 
and let $(P,\mu)$ be the induced symplectic manifold of dimension $2n+2$
constructed via the contact quadruple $(M, N, \alpha, \pi)$.\\
Let $X$ be a  hamiltonian
 vector field on
$M$; i. e.
$$
\L_X\omega=0 \qquad\qquad i(X)\omega =df_X.
$$
Consider the horizontal lift $\bar X$  of $X$ to $N$ defined by
$$
\alpha(\bar X)=0\qquad\qquad \pi_*(\bar X)=X,
 $$
and the lift $\bar{\bar X}$  of $\bar X$ to $P$ defined by
 $$
p_{1*}\bar{\bar X}=\bar X\qquad\qquad ds(\bar{\bar X})=0.
$$
Let $Z$ be the Reeb vector field on $(N,\alpha)$ and let $E$
 be its lift to $P$ defined by
 $$
  p_{1*}E=Z\qquad\qquad ds(E)=0.
$$
\begin{definition}\label{def1.2}
Define the lift  $\tilde X$ of a hamiltonian vector field $X$ on $(M,\omega)$ 
 as the vector field on $P$  defined by:
$$
\tilde X=\bar{\bar X}-(p_1^*\pi^*f_X)\cdot E 
=:\bar{\bar X}-{\widetilde{f_X}}E.
$$
 \end{definition}

\begin{lemma}\label{lem1.4}
 The vector field $\tilde X$ is a hamiltonian vector field on  $(P,\mu)$.
 Furthermore if ${\g}$ is a Lie algebra of vector  fields $X$ on $M$ having a strongly
hamiltonian action, then the set  of vector fields $\tilde X$ on $P$
form an algebra isomorphic to  ${\g}$ and its action on
$(P,\mu)$ is strongly hamiltonian.
\end{lemma}
\begin{proof}
{$i(\tilde X)\mu=i(\bar{\bar X}-{\widetilde{f_X}} E)({\rm
e}^{2s}(p_1^*\pi^*\omega+2ds\wedge p_1^*\alpha))={\rm
e}^{2s}(d\tilde f_X+2ds\tilde f_X)=d({\rm e}^{2s}\tilde f_X)$ which
 shows that $\tilde X$ is hamiltonian and that the hamiltonian function is
 $f_{\tilde X}={\rm e}^{2s}\tilde f_X.$}\\
 Also if $X, Y\in {\g}$:
 \begin{eqnarray*}
   [\tilde X,\tilde Y]&=&[\bar{\bar X}-{\widetilde{f_X}} E,
   \bar{\bar Y}-{\widetilde{f_Y}} E] \\
   &=& \overline{\overline{[X,Y]}}-(\pi\circ p_1)^*\omega(X,Y)E   
   -\widetilde{(Xf_Y)}E+\widetilde{Yf_X} E \\
&=& \overline{\overline{[X,Y]}}- (\pi\circ p_1)^*f_{[X,Y]}E\\
&=& \widetilde{[X,Y]}
\end{eqnarray*}
and 
$$
\{f_{\tilde X},f_{\tilde Y}\}=(\bar{\bar X}-{\widetilde{f_X}}E)({\rm
 e}^{2s}{\widetilde{f_Y}})={\rm e}^{2s}\widetilde{Xf_Y}={\rm e}^{2s}\tilde
 f_{[X,Y]}=f_{\widetilde{[X,Y]}}.
 $$
 \end{proof}
If $C$ is a conformal vector field on $(M,\omega)$ we may assume
$$
{ L}_C\omega =\omega\qquad\qquad di(C)\omega=\omega.
$$
By analogy of what we just did, define the lift $\tilde C_1$ of $C$
 to $(P,\mu)$ by:
 $$
ds(\tilde C_1)=0\qquad p_{1*}\tilde C_1=\bar C+bZ
$$
(i.e. $\pi_*p_{1*}\tilde C_1=C$ and ${\tilde C_1}=\bar{\bar C}+p_1^*b E$). Then
\begin{eqnarray*}
   {L}_{\tilde C_1}\mu &=& d i(\tilde
C_1){\rm e}^{2s}(p_1^*\pi^*\omega+2ds
   \wedge p_1^*\alpha)=d[{\rm e}^{2s}(p_1^* \pi^* i(C)\omega-2 p_1^*b ds)]\\
   &=& {\rm e}^{2s}[p_1^*\pi^*\omega+2ds\wedge p_1^*\pi^*i(C)\omega-2p_1^*db\wedge {ds}] \\
   &=& {\rm e}^{2s}[p_1^*\pi^*\omega+2ds\wedge(p_1^*\pi^*(i(C)\omega)+p_1^*db)].
\end{eqnarray*}
Thus $\tilde C_1$ is a conformal vector field provided:
$$
p_1^*\pi^*i(C)\omega+p_1^*db=p_1^*\alpha.
$$
Or equivalently
$$
 \alpha-\pi^*i(C)\omega=db.
$$
The left hand side is a closed 1-form. If this form is exact 
we are able to lift $C$ to a conformal vector field $\tilde
 C_1$ on $P$. Notice that the rate of variation of $b$ along the
flow  of the Reeb vector field is prescribed:
$$
  Zb=1.
  $$
  A variation of this construction reads as
follows. Let
$$
\tilde C_2=\bar{\bar C}+a E+l\partial_s.
$$
Then:
\begin{eqnarray*}
   { L}_{\tilde C_2}\mu &=& d(i(\bar{\bar
C}+a E+l\partial_s)){\rm e}^{2s}
(p_1^*\pi^*\omega+2ds\wedge p_1^*\alpha) \\
    &=& d\left({\rm e}^{2s}(p_1^*\pi^*(i(C)\omega)-2ads+2lp_1^*\alpha)\right) \\
   &=& {\rm e}^{2s}[p_1^*\pi^*\omega-2da\wedge
ds+2lp_1^*\pi^*\omega+2ds\wedge p_1^*\pi^*i(C)\omega+
   2lds\wedge p_1^*\alpha+2dl\wedge p_1^*\alpha] \\
   &=& {\rm e}^{2s}[(1+2l)p_1^*\pi^*\omega+2ds \wedge
(da+p_1^*\pi^*i(C)\omega+2lp_1^*\alpha)+2dl\wedge p_1^*\alpha].
\end{eqnarray*}
If we choose $l=-1/2$
$$
 { L}_{\tilde C_2}\mu=2{\rm
e}^{2s}ds\wedge(p_1^*\pi^*i(C)\omega-p_1^*\alpha+da).
$$
Thus $\tilde C_2$ is a symplectic vector field on $(P,\mu)$ if the
closed 1-form $\pi^*i(C)\omega-\alpha$ is exact. If this is the case
 the lift $\tilde C_2$ is hamiltonian and
$$
   f_{\tilde C_2}= -a{\rm e}^{2s}.
$$
\begin{lemma}\label{lem1.5}
 If $C$ is a conformal vector field
on $(M,\omega)$, it admits a lift $\tilde C_1$ (resp. $\tilde C_2$)
to $(P,\mu)$ which is conformal  (resp. hamiltonian) if the closed
1-form $\pi^*i(C)\omega-\alpha$  is exact.
\end{lemma}
Let ${\g}$ be an algebra of conformal vector fields on
$(M,\omega)$. Let $X\in{\g}$ be such that
${L}_{X^*}\omega=\omega$ (where $X^*_x=\dfrac{d}{dt}\exp -tX.x|_0$; $x\in
M$). Then ${\g}=\mathbb RX\oplus{\g_1}$, where the vector
fields associated to the elements of $\g_1$, are symplectic.
We
 shall assume here that they are hamiltonian; i.e. $\forall Y\in\g_1$,
 $i(Y^*)\omega=df_Y$.\\
 Consider the lifts of these vector fields to $(P,\mu)$.
\begin{eqnarray*}
 [\tilde X_1^*,\tilde Y^*] &=&[\bar{\bar
X}^*+p_1^*bE,\bar{\bar Y}^*
   -\tilde f_Y E]  \\
    &=& [\bar{\bar X}^*,\bar{\bar Y}^*]-p_1^*\pi^*(Xf_Y)E-p_1^*(\bar Y^*b)E
    +\tilde f_Yp_1^*(Z b) E \\
    &=&
\overline{\overline{[X,Y]}}^{\,*}+[-p_1^*\pi^*\omega(X,Y)
-p_1^*\pi^*\omega(Y,X)+p_1^*\pi^*\omega(X,Y)
    +\tilde f_Y]E \\
&=&\overline{\overline{[X,Y]}}^{\,*}+p_1^*\pi^*(\omega(X,Y)+f_Y)E;\\
    i([X^*,Y^*])\omega&=&-L_{Y^*}i(X^*)\omega=-(i(Y^*)d+di(Y^*))i(X^*)\omega  \\
   &=& -i(Y^*)\omega-d\omega(X,Y)=-d(\omega(X,Y)+f_Y).
\end{eqnarray*}
Hence
$$
[\tilde X_1^*,\tilde Y^*]=\widetilde{[X^*,Y^*]}.
$$
A similar calculation shows that
$$
[\tilde X_2^*,\tilde Y^*]=\widetilde{[X^*,Y^*]}.
$$
Notice as before that $
L_{E}\mu=0$ and $L_{\partial_s}\mu=-2\mu$.
\begin{proposition}
Let $(M,\omega)$ be the first term of a contact quadruple $(M,N,\alpha,\pi)$
and let $(P,\mu)$ be the associated induced symplectic manifold. Then
\begin{enumerate}
\item[(i)] If $G$ is a connected Lie group
acting in a strongly  hamiltonian way on $(M,\omega)$, this action
lifts to a strongly  hamiltonian action of $\tilde G$ (= universal
cover of $G$) on  $(P,\mu)$.
\item[(ii)] If $X$ is a conformal
vector field on $(M,\omega)$ it  admits a conformal (resp.
symplectic) lift to $(P,\mu)$ if the closed  1-form
$\pi^*(i(X)\omega)-\alpha$ is exact. The symplectic lift is  in fact
hamiltonian.
\item[(iii)] The vector field $E$ on $P$ is
hamiltonian and  the vector field $\partial_s$ is conformal.
\end{enumerate}
\end{proposition}
\begin{corollary}\label{corol1.3}
If $(M,\omega)$ admits a transitive hamiltonian action $(P,\mu)$
admits a transitive conformal action. If $(M,\omega)$ admits a
transitive conformal (hamiltonian) action then so does $(P,\mu)$.
\end{corollary}
The stability of the class of conformally homogeneous spaces under
this construction leads us to the study of these spaces.

\section{Conformally homogeneous symplectic manifolds}\label{sec:3}

 \begin{definition}\label{def1.3}
  Let $(M,\omega)$ be a smooth
connected $2n\geq 4$ dimensional  symplectic manifold. A connected
Lie group $G$ is said to  {\bf{act conformally on $(M,\omega)$}} if 
\begin{enumerate}
\item[(i)]$\forall g\in G$,  $g^*\omega=c(g)\omega$ 
\item[(ii)]
There exists at least
one $g\in G$ such that  $c(g)\neq 1$.
\end{enumerate}
\end{definition}
As $\omega$ is closed $c(g)\in \R$; also $c: G\rightarrow \R$ is a
character of $G$. Let $G_1=\ker c$; it is a closed, normal,
 codimension 1 subgroup of $G$.\\
  Let $\g$ (resp. $\g_1$) be the Lie algebra of $G$ (resp. $G_1$).
  Then there exists $0\neq X\in\g$ such that
   $$
\g=\g_1\oplus \R X\qquad \mbox{ (and) }\qquad c_*(X)=1.
$$
The 1-parametric group $\exp tX$ is such that
$$
(\exp tX)^*\omega={\rm e}^t\omega
$$
 and  this group $\exp tX$ is thus isomorphic to $\R$. Hence the group   $G_1$ is connected
and if $G$ is simply connected so is $G_1$. If $X^*$ is the fundamental
vector field on $M$ associated to $X$, remark that $L_{X^*}\omega=-\omega$
since $X^*_x=\frac{d}{dt}\exp -tX\cdot x\vert_0$.
\begin{definition}\label{def1.4}
A symplectic manifold $(M,\omega)$ of dimension $2n\geq 4$ is called
{\bf{conformal homogeneous}} if there exists a Lie group $G$ acting
conformally and transitively on $(M,\omega)$.
\end{definition}
We assume $M$ and $G$ connected. Then $\tilde G$ (= the universal
cover of $G$) is the semi direct product of $\tilde G_1$ (= the
universal cover of $G_1$) by $\R$.\\
By transitivity the orbits of $G_1$ are of dimension $\geq 2n-1$. So
there are two cases
\begin{enumerate}
\item[(i)] The maximum  of the dimension of the $G_1$ orbits is
$(2n-1)$
 \item[(ii)] $G_1$ admits an open orbit.
\end{enumerate}
{\bf Case (i)} By transitivity the dimension of
all $G_1$ orbits is   $(2n-1)$. If we write as above $\g=\g_1\oplus \R X$,
the vector field $X^*$   is everywhere transversal to the $G_1$
orbits. In particular it is  everywhere $\neq 0$. Since $\g_1$ is an
ideal in $\g$, the group   $\exp tX$ permutes the $G_1$ orbits.
Clearly if $\theta_1$ is a
  $G_1$ orbit, $\bigcup_{t\in \R}\exp tX\cdot\theta_1=M$.\\
The restriction $\omega|_{T_x\theta_1}$ has rank $(2n-1)$. Let
 $Z_x$ span the radical of $\omega|_{T_x\theta_1}$ and let
$\alpha:=-i(X^*)\omega\neq 0$ (so $d\alpha=\omega$). As $\alpha_x(Z_x)\neq 0$, we normalize
$Z_x$ so that $\alpha_x(Z_x)=1$. Then
$$
T_xM=\R X^*\oplus
T_x\theta_1=\R X^*\oplus(\R Z_x\oplus\ker \underline{\alpha}_x)
$$
if $\underline{\alpha}_x=\alpha_x|_{T_x\theta_1}$. If $j:
\theta_1\rightarrow M$ denotes the canonical injection
$$
\underline{\alpha}\wedge(d\underline{\alpha})^{n-1}=j^*(\alpha\wedge(\omega)^{n-1})\neq
0.
$$
Thus the orbit $\theta_1$ is a contact manifold, and $Z$ is the
Reeb vector field.\\
Notice that
\begin{eqnarray*}
(L_{X^*}\alpha)(X^*) &=& X^*\alpha(X^*)=0 \\
(L_{X^*}\alpha)(Y^*) &=& X^*\alpha(Y^*)-\alpha([X^*,Y^*]) \\
&=& -X^*\omega(X^*,Y^*)+\omega(X^*,[X^*,Y^*]) \\
&=& -(L_{X^*}\omega)(X^*,Y^*)=\omega(X^*,Y^*)=-\alpha(Y^*)
\end{eqnarray*}
for any $Y\in \g_1$. Hence
$$
L_{X^*}\alpha=-\alpha.
$$
 This says that the various orbits of
$G_1$ have ``conformally''  equivalent contact structure; i. e.
$$
 \underline{\alpha}_{\exp t X\cdot x}(\exp t X_*\cdot
Y^*)={\rm  e}^t\underline{\alpha}_x(Y^*)\qquad Y\in \G_1.
$$
Furthermore
$$
\omega([X^*,Z],Y^*)=X^*\omega(Z,Y^*)-(L_{X^*}\omega)(Z,Y^*)-\omega(Z,[X^*,Y^*])=0
$$
 as $[X^*,Y^*]$is tangent to the orbit. This says that
$[X^*,Z]$ is proportional to $Z$; also
$$
\alpha([X^*,Z])=X^*\alpha(Z)-(L_{X^*}\alpha)(Z)=\alpha(Z)=1
$$
hence $[X^*,Z]=Z$ and thus
$$
(\exp tX)_*Z_x={\rm
e}^tZ_{\exp tX\cdot x}
$$
Finally
\begin{eqnarray*}
\alpha([Y^*,Z]) &=& -\omega(X^*,[Y^*,Z]) \\
&=& -Y^*\omega(X^*,Z)+L_{Y^*}\omega(X^*,Z)+\omega([Y^*,X^*],Z)=0 \\
\omega([Y^*,Z],Y'^*) &=&
Y^*\omega(Z,Y'^*)-L_{Y^*}\omega(Z,Y'^*)-\omega(Z,[Y^*,Y'^*])=0.
\end{eqnarray*}
Hence $[Y^*,Z]$ must be proportional to $Z$ and
thus
$$
 [Y^*,Z]=0
$$
 which says that the Reeb vector is $G_1$ stable.\\
 {\bf Case (ii)} $G_1$ admits an open orbit. We shall assume that
 this orbit coincides with $M$. Thus $(M,\omega)$ is a $G_1$
homogeneous sympletic manifold and $\omega$ is exact.
$$
\omega=d\eta \quad{\rm{where}}\quad \eta:=-i(X^*)\omega.
$$
Assume that the action of $G_1$
is strongly hamiltonian; i. e.  $\forall Y\in \g_1$
\begin{eqnarray*}
i(Y^*)\omega &=& df_Y \\
\{f_Y,f_{Y'}\}&=&-\omega(Y^*,Y'^*)=f_{[Y,Y']}
\end{eqnarray*}
 where $U^*$ denotes the fundamental vector
field associated to  $U\in\g_1$ on $\theta_1$.
Then
\begin{eqnarray*}
L_{Y^*}\eta&=&-L_{Y^*}i(X^*)\omega=-(L_{Y^*}i(X^*)-i(X^*)L_{Y^*})\omega
=-i([Y^*,X^*])\omega\\
&=&df_{[X,Y]}=df_{DY}
\end{eqnarray*}
if $D={\rm ad}X|_{\g_1}.$\\
We also have$L_{X^*}\eta=-\eta$.

 By Kostant's theorem we
may identify $M$ (up to a covering) with a
 coadjoint orbit $\theta_1$ of $G_1$.\\
 Let $\xi\in\theta_1$, let
$\pi:G_1\rightarrow\theta_1:g_1\rightarrow  g_1\cdot\xi={\rm
Ad}^*g_1\xi$  and let $H_1$ be the stabilizer of  $\xi$ in $G_1$. It
is no restriction to assume $X^*_\xi=0$ (since one can replace $X$ by $X+Y$
for any $Y\in \g_1$ and any tangent vector at $\xi$ can be written in 
the form $Y^*_\xi$).\\
 Assuming  $G$ (hence $G_1$) to
be connected and simply connected the derivation  $D$ exponentiates to
a 1-parametric automorphism group of $\g_1$ given by ${\rm e}^{tD}$ and these
``exponentiate'' to a 1-parametric automorphism group
of $G_1$ which will be denoted $a(t)$. The product law in
$G=G_1\cdot \R$ reads:
$$
(g_1,t_1)(g_2,t_2)=(g_1a(t_1)g_2,t_1+t_2).
$$
 As $X^*_\xi=0$ we have:
 \begin{eqnarray*}
(1,t)\cdot\xi &=& \xi \\
(1,t)(g_1,0)\cdot\xi&=& (a(t)g_1,t)\xi=(a(t)g_1,0)(1,t)\cdot\xi\\
&=&(a(t)g_1,0)\cdot\xi=(a(t)g_1\circ g_1^{-1},0)(g_1,0)\cdot\xi.
\end{eqnarray*} In particular if
$g_1\in H_1$ (= stabilizer of $\xi$ in $G_1$) $a(t)g_1\in H_1$;
hence if $Y\in \h_1$ (= Lie algebra of
$H_1$), $[Y,X]\in \h_1$.\\
Furthermore
\begin{eqnarray*}
(L_{X^*}\omega)(Y_1^*,Y_2^*)&=&-\omega(Y_1^*,Y_2^*)\\
&=&X^*\omega(Y_1^*,Y_2^*)-\omega([X^*,Y_1^*],Y_2^*)-\omega(Y_1^*,[X^*,Y_2^*]).
\end{eqnarray*}
 The above relation at
$\xi$ reads:
$$
\omega_\xi(Y_1^*,Y_2^*)=\omega_\xi([X,Y_1^*],Y_2^*)+\omega_\xi(Y_1^*,[X,Y_2^*]).
$$
 But on $\theta_1$, $\omega$ is the Kostant-Souriau symplectic
form;  hence
 \begin{eqnarray*}
\langle\xi,[Y_1,Y_2]\rangle &=& \langle\xi,D[Y_1,Y_2]\rangle \\
\langle\xi-\xi\circ D,[Y_1,Y_2]\rangle &=& 0
\end{eqnarray*}
That is $\xi-\xi D$ vanishes identically on the
derived algebra
$\g'_1$.\\
Conversely suppose we are given an algebra $\g_1$, an element
$\xi\in \g_1^*$ and a derivation $D$ of $\g_1$ such that
$$
 \xi-\xi\circ D\qquad \mbox{ vanishes on }\g'_1.
$$
Then, if, as above, $H_1$ denotes the stabilizer of $\xi$ in $G_1$
and $\h_1$ its Lie algebra, one observes that $Y\in \h_1$ implies 
$DY\in \h_1$. \\
On the orbit $\theta_1=G_1\cdot\xi=G_1/H_1$ define the vector field
${\hat{X}}$ at $\tilde\xi=g_1\cdot\xi$ by:
$$
{\hat{X}}_{\tilde\xi}=\frac{d}{dt}a(-t)g_1\cdot
g_1^{-1}\cdot\tilde\xi|_{t=0}.
$$
 This can be expressed in a
nicer way as:
$$
\langle {\hat{X}}_{\tilde\xi=g_1\xi},Z\rangle=\frac{d}{dt}\langle
a_{-t}(g_1)g_1^{-1}g_1\xi,Z\rangle|_0=\frac{d}{dt}\langle
a_{-t}(g_1)\xi,Z\rangle|_0
$$
for $Z\in \g_1$
\begin{eqnarray*}
{\rm Ad}\,a_{-t}(g_1^{-1})Z &=& \frac{d}{ds}
a_{-t}(g_1^{-1}){\rm e}^{sZ}a_{-t}(g_1)|_0
=\frac{d}{ds}a_{-t}(g_1^{-1} a_t{\rm e}^{sZ} g_1)|_0 \\
&=& \frac{d}{ds}a_{-t}(g_1^{-1}{\rm e}^{s{\rm
e}^{tD_Z}}g_1)|_0=a_{-t*}{\rm Ad}\,g_1^{-1}{\rm e}^{tD}Z      \\
\frac{d}{dt}{\rm Ad}\,a_{-t}(g_1^{-1})Z|_0 &=& -D\circ {\rm
Ad}g_1^{-1}Z+{\rm Ad}\,g_1^{-1}DZ
\end{eqnarray*}
i. e.
$$
{\hat{X}}_{\tilde\xi=g\xi}=-\xi\circ D\circ {\rm
Ad}\,g_1^{-1}+\xi\circ {\rm Ad}\,g_1^{-1}\circ D.
$$
Observe that this
expression has a meaning; indeed if we assume that  $g\in H_1$ (=
stabilizer of $\xi$)
$$
 \xi\circ {\rm Ad}\,g^{-1}=\xi.
$$
Also if $Y\in\h_1,\quad
 \frac{d}{ds}\langle\xi\circ D\circ {\rm
Ad}{\rm e}^{-sY}\rangle|_s =
-\langle\xi\circ D\circ adY\circ{\rm Ad}{\rm e}^{-sY},Z\rangle=0$
so that $
\langle\xi\circ D\circ{\rm Ad}{\rm
e}^{-sY},Z\rangle=\langle\xi\circ D,Z\rangle$.\\
Thus ${\hat{X}}_\xi=0$ and, if $h\in H_1$:
\begin{eqnarray*}
{\hat{X}}_{\tilde\xi=g\cdot\xi=g\cdot h\cdot\xi}&=&\xi\circ D\circ{\rm Ad}\,
h^{-1}\circ{\rm Ad}\,y^{-1}+\xi\circ{\rm Ad}\,h^{-1}\circ{\rm Ad}\,g^{-1}D \\
&=& \xi\circ D\circ{\rm Ad}\,g^{-1}+\xi\circ{\rm Ad}\,g^{-1}D.
\end{eqnarray*}
Furthermore if $Y\in\g_1$:
$$
[Y^*,{\hat{X}}]_{\tilde\xi}=(L_{Y^*}{\hat{X}})_{\tilde\xi}=\frac{d}{dt}
(\varphi^{Y^*}_{-t*}{\hat{X}}_{\varphi^{Y^*}_t\tilde\xi})|_0=
 -(DY)_{\tilde\xi}^*.
$$
Hence, if $Y_1,Y_2\in \g_1$:
\begin{eqnarray*}
(L_{\hat{X}}\omega)_\xi(Y_1^*,Y_2^*)&=&
{\hat{X}}_\xi\omega(Y_1^*,Y_2^*)-\omega((DY_1)^*,Y_2^*)-
\omega(Y_1^*,(DY_2)^*)\\
&=& \langle-\xi,
D[Y_1,Y_2]\rangle=-\langle\xi,[Y_1,Y_2]\rangle=-\omega_\xi(Y_1^*,Y_2^*)
 \end{eqnarray*}
 and similarly at any other point, so that ${\hat{X}}$ is a conformal
 vector field ($L_{\hat{X}}\omega=-\omega$). We conclude by
\begin{proposition}  Let $(M,\omega)$be a smooth
connected $2n (\geq 4)$ dimensional  symplectic manifold which is
conformal homogeneous and let $G$  denote the connected component of
the conformal group. Then
\begin{enumerate}
\item[(i)]$G$
admits a codimension 1 closed, connected, invariant  subgroup $G_1$
which acts symplectically on $M$ and $G/G_1=\R$.
\item[(ii)]If the maximum dimension of the $G_1$ orbits is $(2n-1)$
$M$ is a union of $(2n-1)$ dimensional $G_1$ orbits; each of these
orbits is a contact manifold.
\item[(iii)] If $G_1$ acts
transitively on $M$ in a strongly  hamiltonian way, $M$ is a covering
of a $G_1$ orbit $\theta$ in $\g_1^*$ (= dual of  the Lie
algebra $\g_1$ of $G_1$). Furthermore if $\xi\in\theta$, there exists a
derivation $D$ of $\g_1$ such that
$$
 \xi-\xi\circ D
$$
 vanishes on the derived algebra. Conversely if we are given
an  element $\xi\in\g_1^*$ and a derivation such that $\xi-\xi\circ
D$  vanishes on the derived algebra, the orbit $\theta$ has the
structure of a conformal homogeneous symplectic manifold.
\end{enumerate}
\end{proposition}

\section{Induced connections}\label{sec:4}

We consider the situation where we have  a smooth symplectic manifold
$(M,\omega)$   of dim $2n$, a contact quadruple $(M,N,\alpha,\pi)$
and the corresponding induced symplectic manifold $(P,\mu)$.

Let as before $Z$ be the  Reeb vector field on the contact manifold
$(N,\alpha)$ (i.e. $i(Z)d\alpha=0$ and $\alpha(Z)=1$).
At each point $x\in N$,
$\Ker ({\pi_*}_x)=\R Z$ and $L_Z \alpha = 0$.\\
Recall that $P=N\times \R$ and $
\mu=2e^{2s}~ds\wedge p_1^*\alpha +e^{2s}~ dp_1^*\alpha
$
where $s$ is the variable along $\R$ and  
$p_1: P \rightarrow N$ the projection on the first factor.

Let $\nabla$ be a smooth symplectic connection on $(M,\omega)$.
We shall now define a connection $\nabla^P$ on $P$ induced by $\nabla$.

Let us first recall some notations:\\
Denote by $p$ the projection $p=\pi\circ p_1:P\rightarrow M$.\\
If $X$ is a vector field on $M$, ${\bar{\bar{X}}}$
 is the vector field on $P$ such that
$$
(i)~~p_*{\bar{\bar{X}}}=X \quad\quad  (ii)~~(p_1^*\alpha)({\bar{\bar{X}}})=0
\quad\quad  (iii)~~ds({\bar{\bar{X}}})=0.
$$
We denote by $E$ the vector field on $P$ such that
$$
(i)~~p_{1*}E=Z \quad\quad  (ii)~~ds(E)=0.
$$
Clearly  the values at any point of $P$ of the vector fields 
 ${\bar{\bar{X}}},E, S=\partial_s$ span the tangent space to 
$P$ at that point and we have 
$$
[E,\partial_s]=0 \quad
[E,{\bar{\bar{X}}}]=0Ê\quad
[\partial_s,{\bar{\bar{X}}}]=0 \quad
[{\bar{\bar{X}}},{\bar{\bar{Y}}}]={\overline{\overline{[X,Y]}}}-p^*\omega(X,Y)E.
$$
The formulas for $\nabla^P$ are:
\begin{eqnarray*}
\nabla^P_{{\bar{\bar{X}}}} {\bar{\bar{Y}}}&=&\overline{\overline{{\nabla}_X Y}}
   -\half p^*(\omega(X,Y))E- p^*({\hat{s}}(X,Y))\partial_s
      \\[2mm]
\nabla^P_{E} {\bar{\bar{X}}}&=&\nabla^P_{{\bar{\bar{X}}}}E =2
    \overline{\overline{{\sigma}X}} + p^*(\omega(X,u))\partial_s
      \\[2mm]
\nabla^P_{\partial_s}{\bar{\bar{X}}}&=& \nabla^P_{{\bar{\bar{X}}}}{\partial_s}=\overline{\overline{X}}
      \\[2mm]
\nabla^P_E E&=&p^*f \,{\partial_s}-2\overline{\overline{U}} 
      \\[2mm]
\nabla^P_E{\partial_s}&=&\nabla^P_{\partial_s}E=E
      \\[2mm]
\nabla^P_{\partial_s}{\partial_s}&=&{\partial_s} 
\end{eqnarray*}
where  $f$ is a function on $M$, $U$ is a vector field on $M$, 
${\hat{s}}$ is a symmetric $2$-tensor on
$M$, and $\sigma$ is the endomorphism of $TM$ associated to $s$, hence
${\hat{s}}(X,Y)=\omega(X,\sigma Y)$.

Notice first that these formulas have the correct linearity properties
and yield a torsion free linear connection on $P$.
One checks readily that $\nabla^P\mu=0$ so that $\nabla^P$
is a symplectic connection on $(P,\mu)$.

We now compute the curvature $R^{\nabla^P}$ of this connection $\nabla^P$. 
We get
\begin{eqnarray*}
R^{\nabla^P}({\bar{\bar{X}}},{\bar{\bar{Y}}}){\bar{\bar{Z}}}&=&
\overline{\overline{R^{\nabla}(X,Y)Z}}\cr
&&+\overline{\overline{2\omega(X,Y)\sigma Z-\omega(Y,Z)\sigma X
+\omega(X,Z)\sigma Y-{\hat{s}}(Y,Z)X+{\hat{s}}(X,Z)Y}} \cr
&&+p^*[\omega(X,D(\sigma,U)(Y,Z))
-\omega(Y,D(\sigma,U)(X,Z)] \partial_s\cr
&& ~~\cr
R^{\nabla^P}({\bar{\bar{X}}},{\bar{\bar{Y}}})E&=&
\overline{\overline{2D(\sigma,U)(X,Y)-2D(\sigma,U)(Y,X)}}\cr
&&+p^*[\omega(X,\half fY-\nabla_YU-2{\sigma}^2Y)
-\omega(Y,\half fX-\nabla_XU-2{\sigma}^2X)]\partial_s\cr
&& ~~\cr
R^{\nabla^P}({\bar{\bar{X}}},E){\bar{\bar{Y}}}&=&
\overline{\overline{2D(\sigma,U)(X,Y)}}-
p^*[\omega(Y,\half fX-\nabla_XU-2{\sigma}^2X) ] \partial_s\cr
&& ~~\cr
R^{\nabla^P}({\bar{\bar{X}}},E)E&=&
2\overline{\overline{\half fX-\nabla_XU-2{\sigma}^2X}}
+p^*[Xf+4s(X,u)] \partial_s\cr
&& ~~\cr
R^{\nabla^P}({\bar{\bar{X}}},{\bar{\bar{Y}}}){\partial_{s}}&=&0 \quad \quad
R^{\nabla^P}({\bar{\bar{X}}},E){\partial_{s}}=0 \cr
R^{\nabla^P}({\bar{\bar{X}}},{\partial_{s}}){\bar{\bar{Y}}}&=&0 \quad\quad
R^{\nabla^P}({\bar{\bar{X}}},{\partial_{s}})E=0 \quad 
R^{\nabla^P}({\bar{\bar{X}}},{\partial_{s}}){\partial_{s}}=0\cr
R^{\nabla^P}(E,{\partial_{s}}){\bar{\bar{X}}}&=&0\quad \quad
R^{\nabla^P}(E,{\partial_{s}})E=0\quad
R^{\nabla^P}(E,{\partial_{s}}){\partial_{s}}=0
\end{eqnarray*}
where
$$
D(\sigma,U)(Y,Y'):=(\nabla_Y\sigma)Y'+\half \omega(Y',U)Y-\half\omega(Y,Y')U.
$$
The Ricci tensor $r^{\nabla^P}$ of the connection $\nabla^P$ is given by
\begin{eqnarray*}
r^{\nabla^P}({\bar{\bar{X}}},{\bar{\bar{Y}}})&=&
              r^{\nabla}(X,Y)+ 2(n+1) {\hat{s}}(X,Y)\cr
r^{\nabla^P}({\bar{\bar{X}}},E)&=&
            -(2n+1)\omega(X,u) -2\Tr[Y\rightarrow (\nabla_Y\sigma)(X)]\cr
r^{\nabla^P}({\bar{\bar{X}}},{\partial_s})&=&0\cr
r^{\nabla^P}(E,E)&=&4\Tr (\sigma^2)-2nf
            +2\Tr[X\rightarrow \nabla_XU ] \cr
r^{\nabla^P}(E,{\partial_{s}})&=&0\cr
r^{\nabla^P}({\partial_{s}},{\partial_{s}})&=&0
\end{eqnarray*}

\begin{theorem}\label{th:Ricciflat}
In the framework described above, $\nabla^P$ is a symplectic connection
on $(P,\mu)$ for any choice of ${\hat{s}},U$ and $f$.
The vector field $E$ on $P$ is affine 
( $L_{\tilde{E}}\nabla^P=0$) and symplectic ( $L_{\tilde{E}}\mu=0$);
the vector field $\partial_s$ on $P$ is  affine and conformal 
($L_{\partial_s}\mu=2\mu$).

Furthermore, choosing
\begin{eqnarray*}
{\hat{s}}&=&\frac{-1}{2(n+1)}r^\nabla\cr
{\underline{U}}:&=&\omega(U,\cdot)=\frac{2}{2n+1}\Tr[Y\rightarrow\nabla_Y\sigma]\cr
f&=&\frac{1}{2n(n+1)^2}\Tr (\rho^\nabla)^2 +\frac{1}{n}\Tr [X\rightarrow \nabla_XU].
\end{eqnarray*}
we have:
\begin{itemize}
\item  the connection $\nabla^P$ on 
$(P,\mu)$ is Ricci flat (i.e. has zero 
Ricci tensor); 
\item if the symplectic connection $\nabla$ on $(M,\omega)$ is of
Ricci type, then 
the connection $\nabla^P$ on $(P,\mu)$ is flat.
\item if the connection $\nabla^P$ is locally symmetric, the connection
$\nabla$ is of Ricci type, hence  $\nabla^P$ is flat.
\end{itemize}
\end{theorem}
\begin{proof}
The first point is an immediate consequences of the formulas above for 
$r^{\nabla^P}$. The second point is a consequence of the differential identities
satisfied by the Ricci type symplectic connections  
(which appear in M.~Cahen, S.~Gutt, J.~Horowitz and J.~Rawnsley,
Homogeneous symplectic manifolds with Ricci-type curvature,
\textit{J. Geom. Phys.} \textbf{38} (2001) 140--151).

The third point comes from the fact that
$(\nabla^P_{\bar{\bar{Z}}} R^{\nabla^P})
({\bar{\bar{X}}},{\bar{\bar{Y}}}){\bar{\bar{T}}}$ contains only one term in 
$E$ whose coefficient is $\half W^{\nabla^P}(X,Y,T,Z)$.
\end{proof}

\section{A reduction construction}\label{sec:5}
We present here a  procedure to construct symplectic connections
on some reduced symplectic manifolds; this is a generalisation
of the construction given by P. Baguis and M. Cahen [
 Lett.~Math.~Phys.~57 (2001), pp. 149-160].\\
 Let $(P,\mu)$ be a symplectic manifold of dimension $(2n+2)$.
 Assume $P$ admits a  complete conformal vector field $S$:
$$
L_S\mu=2\mu;\qquad{\rm{define~~}}\alpha:=\half i(S)\mu\qquad {\rm{so~that~}}
d\alpha=\mu.
$$
Assume also that $P$ admits a   symplectic vector field $E$
commuting with $S$
$$
L_E\mu=0\qquad\qquad [S,E]=0 \qquad\qquad (\Rightarrow L_E\alpha=0).
$$
Then $S\mu(S,E)=(L_S\mu)(S,E)=2\mu(S,E)$, so
 if $x$ is a point of $P$ where $\mu_x(S,E)\neq 0$ and if $s$ is a
 parameter along the integral line $\gamma$ of $S$ passing through
$x$ and taking value $0$ at $x$, we have $
\mu_{\gamma(s)}(S,E)={\rm e}^{2s}\mu_x(S,E).$\\
Assume $P':=\{x\in P|\mu_x(S,E)>0\}\neq\emptyset$ and let:
$$
\Sigma=\{x\in P\,|\,\mu_x(S,E)=1\}=\{x\in P\,|\,f_E(x)=\half\}
$$
where $f_E= - i(E)\alpha=-\half \mu(S,E)$ so that $df_E=- L_E\alpha+ i(E)d\alpha=
i(E)\mu$.\\
Thus $\Sigma\neq\emptyset$ and it is a closed hypersurface (called
the constraint hypersurface). Remark that $P'\cong \Sigma\times \R$.\\
 The tangent space to the hypersurface $\Sigma$ is given
 by 
 $$
 T_x\Sigma=\ker (df_E)_x= \ker (i(E)\mu)_x=E^{\perp_\mu}.
 $$
 The restriction of
$\mu_x$ to $T_x\Sigma$ has rank $2n-2$ and a radical spanned by
$E_x$.\\
Remark thus that the restriction of $\alpha$ to $\Sigma$ is a contact
$1$--form on $\Sigma$.

Let $\sim$ be the equivalence relation defined on $\Sigma$ by the
flow of $E$. Assume that the quotient $\Sigma/\sim$
has a $2n$ dimensional manifold $M$ structure 
so that $\pi:\Sigma\rightarrow\Sigma/\sim=M$
is a smooth submersion.\\
Define on $\Sigma$ a ``horizontal'' distribution of dimension $2n$,
$\mathcal H$, by
$$
\mathcal H=>E,S<^{\perp_\mu},
$$
and remark that $\pi_{*\vert_{{\mathcal{H}}_y}}:{{\mathcal{H}}_y}\rightarrow T_{x=\pi(y)}M$
is an isomorphism.\\
Define as usual the reduced 2-form $\omega$ on $M$ by
$$
\omega_{x=\pi(y)}(Y_1,Y_2)=\mu_y(\bar Y_1,\bar Y_2)
$$
where $\bar Y_i$ $(i=1,2)$ is defined by (i) $\pi_*\bar Y_i=Y_i$
(ii) $\bar Y_i\in \H_y$.\\
Notice that $\pi_*[E,\bar Y]=0$,
and $\mu(S,[E,\bar Y])=-L_E\mu(S,\bar Y)+E\mu(S,\bar Y)=0$  hence 
$$[E,\bar Y]=0.$$
The definition of $\omega_x$ does not depend on the choice of $y$.
Indeed
$$
E\mu(\bar Y_1,\bar Y_2)=L_E\mu(\bar Y_1,\bar Y_2)+\mu([E,\bar
Y_1],\bar Y_2)+\mu(\bar Y_1,[E,\bar Y_2])=0.
$$
Clearly $\omega$ is of maximal rank $2n$ as $\H$ is a symplectic
subspace. Finally
\begin{eqnarray*}
  \pi^*(d\omega(Y_1,Y_2,Y_3)) &=& \cyclic_{123} ~(Y_1\omega(Y_2,Y_3)-\omega([Y_1,Y_2],Y_3) )\\
   &=& \cyclic_{123} ~ (\bar Y_1\mu(\bar Y_2,\bar
   Y_3)-\mu(\overline{[Y_1,Y_2]}, \overline{Y}_3))
\end{eqnarray*}
and
$$
[\bar Y_1,\bar Y_2]=\overline{[Y_1,Y_2]}+\mu(S,[\bar Y_1,\bar Y_2])E.
$$
Hence $\omega$ is closed and thus symplectic. Clearly 
$\pi^*\omega=\mu_{\vertÑ\Sigma}=d(\alpha_{\vertÑ\Sigma})$.

\begin{remark}  The symplectic manifold $(M,\omega)$ is the first element of
a contact quadruple $(M,\Sigma,\half\alpha_{\vert_\Sigma},\pi)$  and the associated
symplectic $(2n+2)$--dimensional manifold is $(P',\mu_{\vert_{P'}})$.
\end{remark}

We shall now consider the reduction of a connection.
Let $(P,\mu),E,S,\Sigma, M,\omega$ be as above.
Let $\nabla^P$ be a symplectic connection on $P$ and assume that
the vecor field $E$ is affine ($L_E\nabla^P=0$).

Then define a connection $\nabla^{\Sigma}$ on $\Sigma$ by
$$
\nabla_A^{\Sigma}B:=\nabla_A^PB-\mu(\nabla_A^PB,E)S=\nabla_A^PB+\mu(B,\nabla_A^PE)S.
$$
Then:
\begin{eqnarray*}
  \nabla_A^\Sigma B-\nabla_B^\Sigma A-[A,B] &=& (\mu(B,\nabla_A^PE)-\mu(A,\nabla_B^PE))S \\
   &=& (\mu(B,\nabla_E^PA+[A,E])-\mu(A,\nabla_E^PB+[B,E]))S \\
   &=& (E\mu(B,A)-\mu(B,[E,A])-\mu([E,B],A))S \\
   &=& (L_E\mu(B,A))S=0.
\end{eqnarray*}
Also
\begin{eqnarray*}
  (L_E\nabla^\Sigma)_A B&=&[E,\nabla_A^PB+\mu(B,\nabla_A^PE)S]\\
  &&~\quad-\nabla^P_{[E,A]}B-
  \mu(B,\nabla^P_{[E,A]}E)S-\nabla_A^P[E,B]-\mu([E,B],\nabla_A^PE)S\\
  &=&(L_E\nabla^P)_AB+(E\mu(B,\nabla_A^PE)-\mu(B,\nabla^P_{[E,A]}E)-\mu([E,B],\nabla_A^PE))S\\
  &=&(L_E\mu)(B,\nabla_A^PE)+\mu(B,[E,\nabla_A^PE]-\nabla^P_{[E,A]}E)S=0
\end{eqnarray*}
i. e. $\nabla^\Sigma$ is a torsion free connection and $E$ is an
affine vector field for $\nabla^\Sigma$.

 Define a connection
$\nabla^M$ on $M$ by:
$$
\overline{\nabla_{Y_1}^MY_2}(y)=\nabla_{\bar Y_1}^\Sigma\bar
Y_2(y)-\mu(\bar Y_2,\nabla_{\bar Y_1}^P S)E.
$$
If $x\in M$, this definition does not depend on the choice of $y\in
\pi^{-1}(x)$. Also
\begin{eqnarray*}
  \overline{\nabla_{Y_1}^MY_2}- \overline{\nabla_{Y_2}^MY_1}-\overline{[Y_1,Y_2]}&=&
  \nabla_{\bar Y_1}^\Sigma\bar Y_2- \nabla_{\bar Y_2}^\Sigma\bar
  Y_1-\overline{[Y_1,Y_2]} \\
  &&~\qquad +
  (-\mu(\bar Y_2,\nabla_{\bar Y_1}^PS)+\mu(\bar Y_1,\nabla_{\bar Y_2}^PS))E\\
  &=& \mu(S,[\bar Y_1,\bar Y_2])E+(\mu(\nabla_{\bar Y_1}^P\bar
  Y_2,S)-\mu(\nabla_{\bar Y_2}^P\bar Y_1,S))E=0
\end{eqnarray*}
Finally
\begin{eqnarray*}
  \pi^*((\nabla_{Y_1}^M\omega)(Y_2, Y_3))&=& \pi^*(Y_1\omega(Y_2,Y_3)-\omega(\nabla_{Y_1}^MY_2,Y_3)
  -\omega(Y_2,\nabla_{Y_1}^MY_3)) \\
   &=& \bar Y_1\mu(\bar Y_2,\bar Y_3)-\mu(\nabla_{\bar Y_1}^P\bar
   Y_2+\mu(\bar Y_2,\nabla_{\bar Y_1}^PE)S-\mu(\bar Y_2,\nabla_{\bar
   Y_1}^PS)E,\bar Y_3)\\
   &-&\mu(\bar Y_2,\nabla_{\bar Y_1}^P\bar
   Y_3+\mu(\bar Y_3,\nabla_{\bar Y_1}^PE)S-\mu(\bar Y_3,\nabla_{\bar
   Y_1}^PS)E)\\
   &=&0
\end{eqnarray*}
i. e. the connection $\nabla^M$ is symplectic.
\begin{lemma}\label{lem1.6}
Let $(P,\mu)$ be a symplectic manifold admitting a symplectic
connection $\nabla^P$, a conformal vector field $S$ which is
complete, a symplectic vector field $E$ which is affine and commutes
with $S$. If the constraint manifold $\Sigma=\{x\in P|\mu_x(S,E)=1\}$
is not empty, and if the reduction of $\Sigma$ is a manifold $M$, this
manifold admits a symplectic structure $\omega$ and a natural reduced
symplectic connection $\nabla^M$.
\end{lemma}
In particular
\begin{theorem}\label{th:indred}
Let $(P,\mu)$ be a symplectic manifold admitting 
 a conformal vector field $S$ ($L_S\mu=2\mu$) which is
complete,  a symplectic vector field $E$ which  commutes
with $S$ and assume that, for any $x\in P,~\mu_x(S,E)>0$.
If the reduction of $\Sigma=\{x\in P\,|\,\mu_x(S,E)=1\}$
by the flow of $E$ has a manifold structure $M$ with $\pi:\Sigma\rightarrow M$
a surjective submersion, then
$M$ admits a reduced symplectic structure $\omega$
and $(P,\mu)$ is obtained by induction from $(M,\omega)$
using the contact quadruple $(M,\Sigma, \half i(S)\mu_{\vert_\Sigma},\pi)$.

In particular $(P,\mu)$ admits a Ricci-flat connection.
 \end{theorem}

Reducing $(P,\mu)$ as above and inducing back we see that theorem
\ref{th:Ricciflat} immediately proves this.

{
}


\begin{thebibliography}{9}



\bibitem{CGS} Michel Cahen, Simone Gutt, Lorenz Schwachh\"ofer:
Construction of Ricci-type connections by reduction and induction,
 preprint math.DG/0310375, in The breadth of symplectic and
  Poisson Geometry, Marsden, J.E. and Ratiu, T.S. (eds),
Progress in Math 232,
Birkhauser, 2004.


\bibitem{CS}  M.~Cahen and L.~Schachh\"ofer, Special symplectic connections,
preprint DG0402221.  

\bibitem{Kostant} B. Kostant, Minimal coadjoint orbits and symplectic induction,
in The breadth of symplectic and
  Poisson Geometry, Marsden, J.E. and Ratiu, T.S. (eds),
Progress in Math 232,
Birkhauser, 2004.

\bibitem{Vaisman}   I.~Vaisman,
Symplectic Curvature Tensors
\textit{Monats. Math.} \textbf{100} (1985) 299--327.

\end{thebibliography}
\end{document}